\documentclass[aps,preprint,a4paper,floatfix]{revtex4}

\usepackage{graphicx,color}
\usepackage{amsmath,amssymb,amsfonts,amsthm,amscd,bm}
\usepackage{booktabs}
\usepackage{cancel}

\def\tilde{\widetilde}

\def\*{\star}
\def\[{\left[}
\def\]{\right]}
\def\({\left(}      

\def\){\right)}

\def\frac#1#2{\dfrac{#1}{#2}}
\def\inv#1{\dfrac{1}{#1}}

\def\2pi{\hbox{$2\pi i$}}

\def\dsl{\raise.15ex\hbox{/}\kern-.57em\partial}
\def\Dsl{\,\raise.15ex\hbox{/}\mkern-.13.5mu D}
\def\Li{{\rm Li}}
\def\2pi{\hbox{$2\pi i$}}

\def\dsl{\raise.15ex\hbox{/}\kern-.57em\partial}
\def\Dsl{\,\raise.15ex\hbox{/}\mkern-.13.5mu D}
\def\barray{\begin{eqnarray}}
\def\earray{\end{eqnarray}}
\def\beq{\begin{equation}}
\def\eeq{\end{equation}}

\def\phi{\Phi}

\def\Li{{\rm Li}}

\def\AA{\leavevmode\setbox0=\hbox{h}
\dimen0=\ht0 \advance\dimen0 by-1ex\rlap{\raise.67\dimen0\hbox{\char'27}}A}

\def\yzero{y_\bullet}

\makeatletter
\def\iddots{\mathinner{\mkern1mu\raise\p@
\vbox{\kern7\p@\hbox{.}}\mkern2mu
\raise4\p@\hbox{.}\mkern2mu\raise7\p@\hbox{.}\mkern1mu}}
\makeatother

\def\fref#1{FIG. \ref{#1}}
\def\tref#1{TABLE \ref{#1}}

\theoremstyle{plain}

\theoremstyle{remark}
\newtheorem{remark}{Remark}

\begin{document}

\title{
Statistical and other properties of Riemann zeros 
based on an  explicit equation for the \boldmath{$n$}-th zero
on the critical line
}
\author{
Guilherme Fran\c ca\footnote{guifranca@gmail.com}  
and  Andr\'e  LeClair\footnote{andre.leclair@gmail.com}
}
\affiliation{Cornell University, Physics Department, Ithaca, NY 14850} 

\begin{abstract}

We show that there are an  infinite  number of Riemann zeros on the critical 
line,  enumerated by the positive integers $n=1,2,\dotsc$,   whose ordinates 
can be obtained as the solution of a new transcendental equation that depends 
only on $n$. Under weak assumptions,  we show that the number of such zeros 
already saturates the counting formula for the numbers of zeros on the 
entire critical strip.  These results thus constitute a concrete proposal 
toward verifying the Riemann hypothesis.  We perform numerical analyses of 
the exact equation, and its asymptotic limit of large ordinate.   
The starting point is an explicit analytical formula for an approximate 
solution to the exact equation in terms of the Lambert $W$ function.  
In this way, we neither have to use Gram points or deal with violations 
of Gram's law. Our numerical approach thus constitutes a novel method to 
compute the zeros.  Employing these numerical solutions, we verify that 
solutions of the asymptotic version  are accurate enough to confirm 
Montgomery's  and Odlyzko's pair correlation conjectures and also to 
reconstruct the prime number counting function. 

\end{abstract}

\maketitle

\section{Introduction}

Riemann's major contribution to number theory was an explicit formula 
for the arithmetic function $\pi (x)$, which counts the number of primes 
less than $x$, in terms of an infinite sum over the non-trivial zeros of 
the $\zeta (z) $ function,  i.e. roots  $\rho$ of the equation $\zeta (z ) =0$ 
on the \emph{critical strip} $0\leq \Re (z)  \leq 1$ \cite{Edwards}.
It was later proven by Hadamard and de la Vall\' ee Poussin that there are 
no zeros on the line $\Re (z) =1$,  which in turn  proved the Prime Number 
Theorem $\pi (x)  \sim  \Li (x)$.  (See section \ref{sec:prime} for a review.)  Hardy proved that there are an infinite number of zeros on the 
\emph{critical line} $\Re (z) = \tfrac{1}{2}$. The \emph{Riemann hypothesis} 
(RH) was his  statement, in 1859,  that all zeros on the critical strip 
have $\Re(\rho)=\tfrac{1}{2}$, although he was unable to prove it. 
Despite strong numerical evidence of its validity,  it remains unproven 
to this day. Many important mathematical results were proven assuming the RH, 
so it is a cornerstone of fundamental mathematics. Some excellent 
introductions to the RH are \cite{Conrey,Sarnak,Bombieri}.

Throughout this paper, the argument of the $\zeta(z) $ function will be 
the complex number $z= x + iy$,  and zeros will be denoted as $\rho$.
We need only consider the positive $y$-axis,  since if $\rho$ is a zero  
so is its complex conjugate.   The infinite zeros along the critical line can 
be numbered as one moves up the $y$-axis, $\rho_n = \tfrac{1}{2} +i y_n$. 
The first few are $y_1 \simeq 14.1347$,  $y_2 \simeq 21.0220$ 
and $y_3 \simeq 25.0108$.   Although at first sight there doesn't appear to be
any regular pattern to these zeros,  
we will demonstrate in this paper that they have a universal description:   
there are in one-to-one correspondence with the zeros of the cosine function.  

Riemann gave an estimate $N(T)$ for the  average number of zeros
on the \emph{entire critical strip} with imaginary part between  
$0$ and $T$. If $T$ does \emph{not} correspond to the ordinate of a zero,
when $T\to \infty$ we have \cite{Edwards, Titchmarsh}
\beq\label{riemann_counting}
N(T) = \dfrac{T}{2\pi}\log\(\dfrac{T}{2\pi e}\) + \dfrac{7}{8}
+O \(\log T \).
\eeq
This formula was later proven by von Mangoldt, but has it never been proven 
to be valid \emph{on the critical line},  as explicitly stated in 
Edward's book \cite{Edwards}.
Denoting the zeros on the critical line by $N_0(T)$, Hardy and
Littlewood showed that $N_0(T) > C \, T$ and Selberg improved this result
stating that $N_0(T) > C \, T\log T$ for very small $C$. 
Then, Levinson \cite{Levinson} demonstrated that 
$N_0(T) \ge C N(T)$ where $C = \tfrac{1}{3}$. The current most 
precise result is due to Conrey \cite{Conrey2} who improved the last result
demonstrating that $C = \tfrac{2}{5}$.
Obviously, if the RH is true then 
$N_0(T)=N(T) \sim \tfrac{T}{2\pi}\log \tfrac{T}{2\pi} - \tfrac{T}{2\pi}$. 
These statements
are described in \cite[Chapter 11]{Edwards} and \cite[Chapter X]{Titchmarsh}.
The formula \eqref{riemann_counting} can be seen as an asymptotic expansion
of an exact formula due to Backlund, who proved the following result
also \emph{on the critical strip} \cite[Chapter 6]{Edwards}:
\beq\label{backlund}
N(T) = \dfrac{1}{\pi}\vartheta(T) + 1 + S(T),
\eeq
where we have the Riemann-Siegel $\vartheta$ function 
(introduced in section \ref{sec:exact}) and
$S(T) = \tfrac{1}{\pi}\arg\zeta(\tfrac{1}{2}+iT)$. 
Using the well known expansion 
$\vartheta(T) = \tfrac{T}{2}\log\(\tfrac{T}{2\pi e}\)-\tfrac{\pi}{8}+
O\(T^{-1}\)$ one recovers \eqref{riemann_counting} from \eqref{backlund}. 

Montgomery's conjecture that the non-trivial zeros 
satisfy the statistics of
the eigenvalues of random hermitian matrices \cite{Montgomery} led
Berry to propose that the zeros are eigenvalues of a chaotic 
hamiltonian \cite{Berry1},   along the lines of the original 
Hilbert-Polya idea.   
Further developments are in \cite{BerryKeating, BerryKeating2, 
Sierra,  Sierra2,  Bhaduri}. These works focus on $N(T)$,
and  carry out the analysis on the 
critical line, i.e. they essentially assume the validity of the RH.    
A number of interesting analytic results were obtained, emphasizing the 
important  role of the function $\arg \zeta\(\tfrac{1}{2}+iT\)$. 
In a related,  but essentially different approach by Connes based on 
adeles, there exists an 
operator playing the role of the hamiltonian,
which has a continuous spectrum, and the Riemann zeros correspond 
to missing spectral lines \cite{Connes}. We mention these interesting 
works because of the role of $N(T)$ in them, however, we will not be 
pursuing these ideas in this work. For interesting connections of the RH
to physics see \cite{Connes2,Schumayer} (and references therein).
 
Riemann's counting formula \eqref{riemann_counting} 
counts zeros very accurately if one takes into account the term 
$\tfrac{1}{\pi}\arg \zeta\(\tfrac{1}{2}+iT\)$.  Thus, 
it is not a smooth function but jumps by one at each zero on the critical 
line. This ``fluctuating term''  is discussed in some detail in 
\cite{Berry1,BerryKeating}. If in some region of the critical 
strip one can show that the counting formula $N(T)$ correctly counts the zeros 
on the \emph{critical line}, then this proves the RH in this region of 
the strip.
Since it has been shown numerically that the first billion or so zeros all 
lie on the critical line \cite{deLune,Gourdon}, one approach to 
establishing the RH is to develop an asymptotic approximation  and show that 
there are no zeros off of the critical line  for sufficiently
large $y$.    Such an analysis was carried out in \cite{RHLeclair} where the
main outcome  was an asymptotic equation for the $n$-th zero on 
the critical line, $\rho_n = \tfrac{1}{2} + i y_n$, where $y_n$ satisfies 
the transcendental equation \eqref{FinalTranscendence} below.  
The way in which this equation is derived shows that these zeros are in 
one-to-one correspondence with the zeros of the cosine function;   
it is in this manner
that the $n$-dependence arises.   
As will be shown in this paper, the numerical solutions to this equation 
unexpectedly accurately  correspond to the already well known values 
for $y_n$ \cite{Odlyzko}, even for the lowest zeros.

More importantly,  since  these  equations for zeros on the critical line are 
enumerated by the integer $n$,  one can use them to obtain the counting of 
such zeros,  which we continue
to denote  as $N_0 (T)$.    Comparing with
Riemann's counting formula \eqref{riemann_counting} for the number 
of zeros on the entire critical strip,  
we  will argue  that $N_0(T) = N(T)$,  
first  asymptotically, then exactly,   based on the exact equation 
\eqref{exact_eq2}.   

Our work presents a novel method to compute the Riemann zeros.  
We first  obtain an explicit formula as an approximate solution for $y_n$, 
in terms of the Lambert $W$ function. 
Starting from this approximation we  obtain 
accurate numerical solutions of \eqref{FinalTranscendence}, which is 
the simplest approximation to \eqref{exact_eq2}.
We show that these numerical solutions are accurate enough to verify  
Montgomery's  and Odlyzko's pair correlation conjectures, and also to 
reconstruct the prime number counting formula. We emphasize that our numerical 
approach does not make use of Gram points nor the Riemann-Siegel $Z$ function,
and  we  believe is actually simpler than the standard methods.   

Let us anticipate a possible  misunderstanding or  criticism due to the 
resemblance between \eqref{FinalTranscendence} and \eqref{riemann_counting}, 
and also between \eqref{exact_eq2} and \eqref{backlund}. 
We stress that our results were derived  directly on the critical line, 
\emph{without assuming the RH}.  Furthermore, \eqref{FinalTranscendence} 
and \eqref{exact_eq2} are not counting formulas. 
Rather, they are \emph{equations} that determine the imaginary
parts $y_n$'s of the Riemann zeros. In other words, 
the $n$-th Riemann zero is the solution of these equations. 
Whereas the simple equation $\zeta (\rho) =0$ has an infinite number of 
solutions, equations \eqref{FinalTranscendence} and \eqref{exact_eq2}  
have  a {\it single} solution for each $n$.
We remind  the reader 
that formulas \eqref{riemann_counting} and \eqref{backlund} were derived
on the \emph{entire critical strip}, moreover, assuming that 
$T$ is \emph{not} the ordinate of a zero. Thus, it is impossible to 
derive \eqref{FinalTranscendence} from \eqref{riemann_counting}, 
nor  \eqref{exact_eq2} from \eqref{backlund}.       
The equations \eqref{FinalTranscendence} 
and \eqref{exact_eq2}  are new equations that are fundamentally 
different in meaning,  and stronger,  than the  known counting formulas.  
We have been unable to find them  in the literature.

We organize our work as follows. Section \ref{sec:main_formula} contains
our main results. More precisely, we derive an exact equation satisfied 
by each individual Riemann zero on the critical line.  The asymptotic limit 
of this equation is the  equation  first proposed in \cite{RHLeclair},  
however we provide  a more rigorous and thorough  analysis.
In section \ref{sec:lambert} we obtain an approximate solution for the
ordinates of the zeros on the critical line, as an explicit formula.
This provides the starting point to compute accurate numerical solutions, 
shown in section \ref{sec:numerical}. In section \ref{sec:gue} we verify the
Montgomery-Odlyzko pair correlation conjecture, based on our numerical
solutions of the asymptotic version.   
Also, in section \ref{sec:prime} we reconstruct the prime number
counting function,  again based on solutions of the asymptotic 
approximation  of the exact equation.   
Section \ref{sec:numerical_exact} presents some numerical
solutions to the exact equation, which proved to be much more robust
under the numerical methods. Finally, in section \ref{sec:conclusion}, we
present our concluding remarks.

\section{An equation for the  Riemann zeros on the critical line}
\label{sec:main_formula}

In this section  we derive the exact equation \eqref{exact_eq2} for the
$n$-th Riemann zero,  which is our main result.     
In the first sub-section  we present its asymptotic version 
\eqref{FinalTranscendence},  first proposed in 
\cite{RHLeclair},  since it involves more familiar functions;  
this first sub-section should be viewed as following trivially from the 
second sub-section.

\subsection{Asymptotic equation}
\label{sec:asymptotic}

Let us start by  defining the function 
\beq
\label{chidef}
\chi(z) \equiv \pi^{-z/2} \Gamma \(z/2\) \zeta (z) .
\eeq
In quantum statistical physics, this function is the free energy of a 
gas of massless bosonic particles in $d$ spatial dimensions when $z=d+1$,  
up to the overall power of the temperature $T^{d+1}$.     
Under a ``modular'' transformation that exchanges one spatial 
coordinate with Euclidean time,  if one analytically  continues $d$,   
physical arguments \cite{AL}  shows that it must have 
the symmetry 
\beq
\label{chisym}
\chi\(z\) = \chi\(1-z\).
\eeq
This is the fundamental, and amazing, functional equation 
satisfied by the $\zeta(z)$ function, which was proven by  Riemann.
  For several different ways of 
proving \eqref{chisym} see \cite{Titchmarsh}. Now
consider Stirling's approximation, 
$\Gamma (z)  = \sqrt{2 \pi}  z^{z - 1/2}e^{-z}\(1+O\(z^{-1}\)\)$,
where $z=x+iy$,  which is valid for large $y$. Under this
condition we also have
\beq
z^z = \exp\( i\(y\log y + \dfrac{\pi x}{2}\) + x \log y - \dfrac{\pi y}{2} 
+ x + O\(y^{-1}\) \).
\eeq
Therefore, using the polar representation $\zeta = |\zeta| e^{i\arg\zeta}$ and
the above expansions,  we can write
$\chi = A \, e^{i\theta}$ where
\begin{align}
A(x,y) &= \sqrt{2\pi } \, \pi^{-x/2} \(  \dfrac{y}{2} \)^{(x-1)/2} 
e^{- \pi y /4} |\zeta (x + i y)|\(1+O\(z^{-1}\)\) ,  
\label{A_assymp} \\
\theta(x, y) &= \dfrac{y}{2} \log \( \dfrac{y}{2 \pi e} \)  + 
\dfrac{\pi}{4}(x-1) + \arg \zeta(x + i y) + O\(y^{-1}\) .
\label{theta_assymp}
\end{align}
The above approximation is very accurate. For $y$ as low as $100$, 
it evaluates $\chi\(\tfrac{1}{2}  + i y\)$ correctly to one part in $10^6$.

Now let $\rho = x+iy$ be a Riemann zero. Then $\arg \zeta(\rho)$ can be 
well-defined by the limit
\beq \label{deltadef} 
\arg\zeta\(\rho\) \equiv  \lim_{\delta \to 0^+} \arg \zeta\(x+\delta+iy\).
\eeq
Note that $0 < \delta \ll 1$. This limit in general is not zero.
For instance, for the first Riemann zero, 
$\arg \zeta\(\tfrac{1}{2} + i y_1\) \simeq 0.1578739$.
On the critical line $z=\tfrac{1}{2}+it$, if $t$ does not
correspond to the imaginary part of a zero, the well known function 
$S(t) = \tfrac{1}{\pi}\arg\zeta\(\tfrac{1}{2}+it\)$, already mentioned
in connection with \eqref{riemann_counting} and \eqref{backlund}, 
is defined by continuous variation along the straight 
lines starting from $2$, then up to $2+it$ and finally to
$\tfrac{1}{2}+it$, where $\arg\zeta(2)=0$. Assuming the RH, the current 
best bound is given by
$|S(t)|\le\(\tfrac{1}{2}+o(1)\)\tfrac{\log t}{\log\log t}$ for 
$t\to \infty$, proven by Goldston and Gonek \cite{Goldston}. 
On a zero, the standard way to define this term is through the limit
$S(\rho) = \tfrac{1}{2} 
\lim_{\epsilon\to0} \( S\(\rho+i\epsilon\)+S\(\rho-i\epsilon\) \)$.
We have checked numerically that for several zeros on the line, our definition
\eqref{deltadef} gives the same answer as this standard approach.

From \eqref{chidef} we have $\(\chi(z)\)^{*} = \chi\(z^*\)$, thus 
$A(x,-y)=A(x,y)$ and $\theta(x,-y)=-\theta(x,y)$. 
Denoting $\chi\(1-z\) \equiv A' \, e^{-i\theta'}$ this implies that
\beq
A'(x,y)=A(1-x,y), \qquad \theta'(x,y)=\theta(1-x,y).
\eeq
From \eqref{chisym} we also have $|\chi(z)| = |\chi(1-z)|$, 
therefore $A(x,y) = A'(x,y)$ for \emph{any} $z$ on the critical strip.

Now let us consider what happens when we approach a zero
$\rho = x + iy$ through a limit. 
From \eqref{chidef} it follows that $\zeta(z)$ and $\chi(z)$ have
the same zeros on the critical strip, so it is enough
to consider the zeros of $\chi(z)$. From \eqref{chisym} we see 
that if $\rho$ is a zero so is $1-\rho$. Then we clearly have
\footnote{The linear combination in \eqref{sumchi}  was chosen 
to be manifestly symmetric under $z\to 1-z$. Had we taken a different 
linear combination in \eqref{sumchi}, 
then $B= e^{i\theta} + b \, e^{-i\theta'}$ for some constant $b$.
Setting the real and imaginary parts of $B$ to zero gives the two equations
$\cos \theta + b\cos \theta' =0$ and $\sin\theta - b  \sin \theta'=0$.    
Summing the squares of these equations one obtains  
$\cos (\theta + \theta')=-(b+1/b)/2$. However, since $b+1/b >1$,  there 
are no solutions except for $b=1$.}
\beq \label{sumchi}
\lim_{\delta\to 0^{+}}\[ \chi(\rho+\delta) + \chi(1-\rho-\delta) \] = 
\lim_{\delta \to 0^{+}} A(x+ \delta,y) B(x+\delta, y) = 0, 
\eeq
where 
\beq
B(x,y) \equiv e^{i \theta(x,y) }  +  e^{-i \theta'(x,y)}.
\eeq
The second equality in \eqref{sumchi} follows from $A=A'$. 
Then, in the limit $\delta \to 0^{+}$, a zero corresponds 
to $A=0$, $B=0$ or both.
They can simultaneously be zero since they are not independent.
If $B=0$ then $A=0$, since $A \propto |\zeta (z) |$. However, the converse 
is not necessarily true.

Since there is more structure in $B$, let us consider $B=0$.
The general solution of this equation  is given by 
$\theta + \theta' = (2n+1)\pi$, which are a family of curves $y(x)$.
However, since $\chi(z)$ is an analytic function, we know that the 
zeros must be isolated  points rather than curves, and this general 
solution must be restricted. Thus, let us choose the particular solution
\beq \label{particular_sol}
\theta = \theta', \qquad \lim_{\delta\to 0^{+}}\cos\theta = 0.
\eeq
On the critical line, the first equation \eqref{particular_sol} is 
already satisfied. Now, in the limit $\delta\to0^+$, the second 
equation implies
$\theta = \(n+\tfrac{1}{2}\)\pi$, for $n=0,\pm 1, \pm 2, \dotsc$, hence
\beq
\label{almost_final}
n = \dfrac{y}{2\pi}\log\(\dfrac{y}{2\pi e}\) -\dfrac{5}{8}
+ \lim_{\delta\to0^{+}}\dfrac{1}{\pi}\arg \zeta\(\tfrac{1}{2}+i y + \delta\).
\eeq
A closer inspection shows that the right hand side of
\eqref{almost_final} has a minimum in the interval $(-2, -1)$, thus $n$ 
is bounded from below, i.e. $n \ge -1$. 
Establishing the \emph{convention} that zeros 
are labeled by positive integers, $\rho_n = \tfrac{1}{2}+i y_n$ where
$n=1,2,\dotsc$, we must replace $n \to n - 2$ in \eqref{almost_final}. 
Therefore, the imaginary
parts of these zeros are determined from the solution
of the transcendental equation
\beq
\label{FinalTranscendence} 
\dfrac{y_n}{2 \pi}  \log \( \dfrac{y_n }{2 \pi e} \)   
+ \lim_{\delta \to 0^{+}}  \dfrac{1}{\pi} 
\arg  \zeta \( \tfrac{1}{2}+ \delta + i y_n  \) = n - \dfrac{11}{8}
\qquad (n=1,2,\dotsc).
\eeq
In short, we have shown that, asymptotically, there are an  infinite 
number of  zeros on the critical line whose ordinates can be determined by 
solving \eqref{FinalTranscendence} for $y_n$.

Note that,  by comparing with the counting function $N(T)$,   the left hand 
side of \eqref{FinalTranscendence} is a monotonic increasing function 
of $y$,  and the leading term is a smooth function. Possible discontinuities 
can only come from 
$\tfrac{1}{\pi}\arg\zeta\(\tfrac{1}{2}+iy\)$, and in fact, it has a jump 
discontinuity by one whenever $y$ corresponds to a zero. 
However,  if $\lim_{\delta\to0^+}\arg\zeta\(\tfrac{1}{2}+\delta+ iy \)$ is 
well defined,
then the left hand side of equation \eqref{FinalTranscendence} is well 
defined for any $y$ and there is a unique solution for every $n$.
Under this assumption, the number of solutions of equation 
\eqref{FinalTranscendence}, up to height $T$, is given by
\beq \label{counting2}
N_0(T) = \dfrac{T}{2 \pi} \log \( \frac{T}{2 \pi e} \) + \frac{7}{8}  + 
\inv{\pi} \arg \zeta \( \tfrac{1}{2} + i T \) + O\(T^{-1}\).
\eeq
This is so because the zeros are already numbered in 
\eqref{FinalTranscendence}, but the left hand side jumps by one at each 
zero, with values $-1/2$ to the left and $+1/2$ to the right of the zero. 
Thus we can replace $n \to N_0 + 1/2$ and $y_n \to T$, such that the 
jumps correspond to integer values. In this way $T$ will not correspond 
to the ordinate of a zero and $\delta$ can be eliminated.

Let us  now recall  the Riemann-von Mangoldt formula \eqref{riemann_counting}
for the number of zeros on the \emph{critical strip}. It is
the same as the number of zeros on the \emph{critical line} that
we have just found \eqref{counting2}, i.e. $N_0(T) = N(T)$.
This means that our particular solution \eqref{particular_sol}, leading
to equation \eqref{FinalTranscendence}, already saturates the counting 
formula on the whole strip and there are no additional zeros from $A=0$ 
in \eqref{sumchi} nor from the general solution 
$\theta + \theta' = (2n+1)\pi$. This strongly suggests that
\eqref{FinalTranscendence} describes all the non-trivial zeros, which are 
 all on the critical line.

\subsection{Exact equation}
\label{sec:exact}

Let us  now reproduce the same analysis discussed previously but without
an asymptotic expansion. The exact versions of
\eqref{A_assymp} and \eqref{theta_assymp} are
\begin{align}
\label{A_exact}
A(x,y) &= \pi^{-x/2} |\Gamma\(\tfrac{1}{2}(x+iy)\)| |\zeta(x+iy)|, \\
\label{theta_exact}
\theta(x,y) &= 
\arg \Gamma\(\tfrac{1}{2}(x+iy)\) -\dfrac{y}{2}\log\pi + \arg \zeta(x+iy),
\end{align}
where again $\chi(z)=Ae^{i\theta}$ and $\chi(1-z)=A'e^{-i\theta}$, with
$A'(x,y)=A(1-x,y)$ and $\theta'(x,y) = \theta(1-x,y)$. 
The zeros on the critical line correspond to the  
particular solution $\theta=\theta'$ and 
$\lim_{\delta\to 0^+}\cos \theta =0$. Thus 
$\lim_{\delta\to0^+}\theta = \(n+\tfrac{1}{2}\)\pi$ 
and replacing $n \to n-2$, the imaginary parts of these zeros must 
satisfy the exact equation
\beq\label{exact_eq}
\arg\Gamma\(\tfrac{1}{4}+\tfrac{i}{2}y_n\) - y_n \log\sqrt{\pi} 
+ \lim_{\delta \to 0^+} \arg\zeta\(\tfrac{1}{2}+iy_n\) = \(n-\tfrac{3}{2}\)\pi.
\eeq
The Riemann-Siegel $\vartheta$ function is defined by
\beq
\label{riemann_siegel}
\vartheta(t)\equiv \arg\Gamma\(\tfrac{1}{4}+\tfrac{i}{2}t\) - t\log \sqrt{\pi},
\eeq
where the argument is defined such that this function is continuous and
$\vartheta(0)=0$.   Therefore, there are infinite zeros in the form
$\rho_n=\tfrac{1}{2}+iy_n$, where $n=1,2,\dotsc$, whose imaginary
parts \emph{exactly} satisfy the following equation:
\beq\label{exact_eq2}
\vartheta(y_n) + 
\lim_{\delta\to 0^{+}}\arg\zeta\(\tfrac{1}{2}+\delta +iy_n\) = 
\(n-\tfrac{3}{2}\)\pi \qquad (n=1,2,\dotsc).
\eeq
Expanding the $\Gamma$-function in \eqref{riemann_siegel} through 
Stirling's formula, one recovers the asymptotic equation 
\eqref{FinalTranscendence}. 

We now argue that \eqref{exact_eq2}  has a unique solution for 
each $n$. Let $g(y)$ be the function defined by its left hand side (with
$y_n \to y$). The function $g(y)$ is monotonically increasing,  and 
the shift by $\delta$ makes $g(y)$ well-defined  between 
the discontinuous jumps of the $\arg\zeta$ term.
The reason that $\delta$ must be taken  positive
is the following. Near a zero $\rho_n$,   
$\zeta (z)  \approx  (z-\rho_n) \zeta' (\rho_n)
= (\delta +i (y-y_n) ) \zeta'(\rho_n)$.    
This gives $\arg \zeta (z) \approx \arctan ((y-y_n)/\delta ) + 
{\rm const.}$.    Thus,  with $\delta >0$,  as one passes through a zero 
from below,  $S(y)$ 
{\it increases} by $1$ as it should based on its role in the counting 
function $N(T)$.
Thus the equation $g(y) = (n -3/2)\pi$ should have a unique solution 
for every $n$. Under this condition it is valid to replace $y_n \to T$ and 
$n \to N_0 + \tfrac{1}{2}$ into \eqref{exact_eq2}, yielding
the number of zeros on the critical line
\beq
N_0(T) = \dfrac{1}{\pi}\vartheta(T) + 
\dfrac{1}{\pi}\arg\zeta\(\tfrac{1}{2}+iT\) + 1.
\eeq
Therefore, comparing with the exact counting formula on the \emph{whole
strip} \eqref{backlund}, we have $N_0(T) = N(T)$ exactly. This indicates
once again that our particular solution, 
leading to equation \eqref{exact_eq2}, captures all the zeros on the strip,
showing that they should all be on the critical line. 
In summary,  if \eqref{exact_eq2}  has a unique solution for each $n$,  
as we have argued,   then this proves the RH.

\subsection{Further remarks}

\begin{figure}
\centering
\begin{minipage}{.5\textwidth}
  \centering
  \includegraphics[width=1\linewidth]{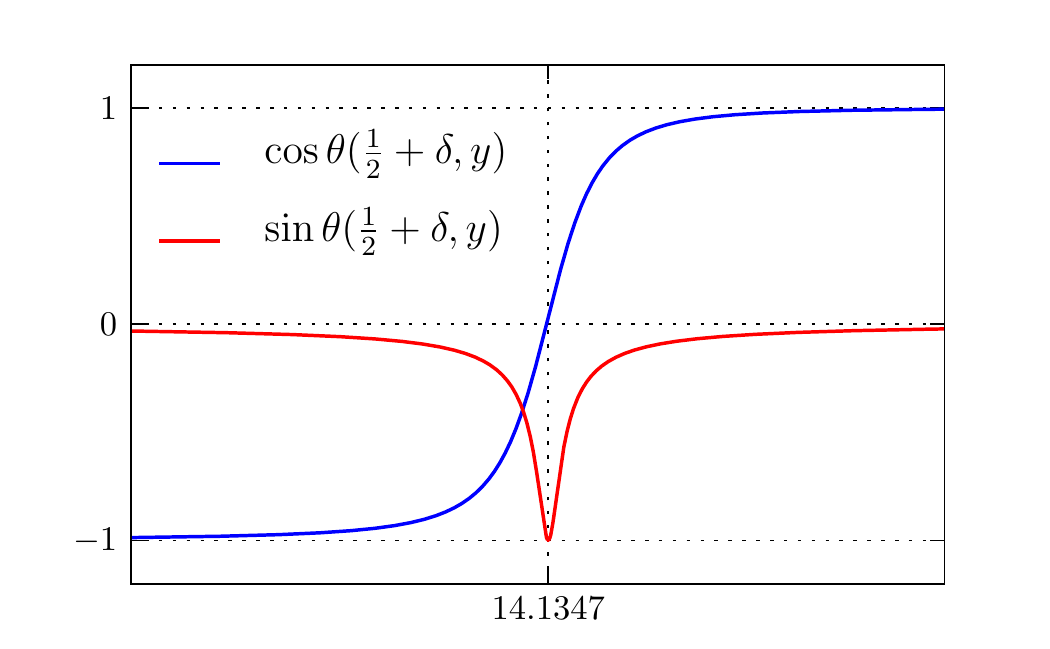}\\[-1em]
(a)
\end{minipage}%
\begin{minipage}{.5\textwidth}
  \centering
  \includegraphics[width=1\linewidth]{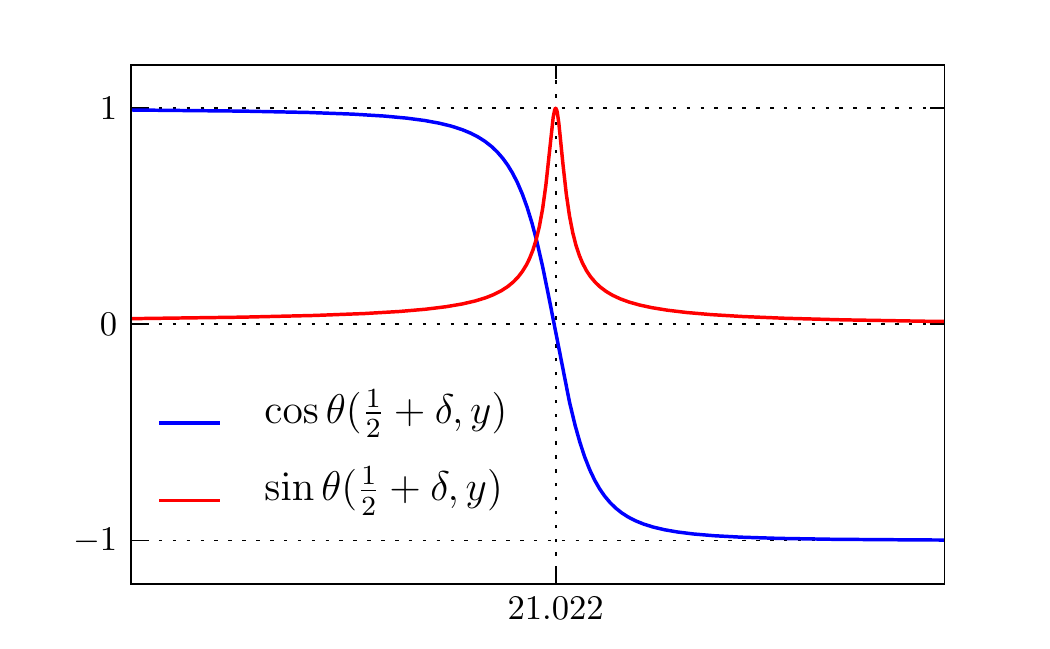}\\[-1em]
(b)
\end{minipage}
\caption{Exactly on a zero we have $\cos\theta = 0$ and
$\sin\theta = \pm 1$. This is illustrated for
the first (a) and second (b) Riemann zeros, respectively. We plot
$\cos\theta\(\tfrac{1}{2}+\delta, y\)$ and 
$\sin\theta\(\tfrac{1}{2}+\delta, y\)$ versus $y$, for
$0 < \delta \ll 1$.}
\label{fig:cos} 
\end{figure}

\begin{remark}\label{simple_zeros}
An important consequence of equation \eqref{exact_eq2}, or its
asymptotic version \eqref{FinalTranscendence}, is
that all of its zeros are simple. This follows from the fact that they 
are in one-to-one correspondence with the zeros of the cosine 
function \eqref{particular_sol}, which are simple.
If the zeros are simple, there is an easier  way to see 
that the zeros correspond to $\cos \theta =0$. On the critical
line $z = \tfrac{1}{2} + i y$, the functional equation 
\eqref{chisym} implies $\chi (z)$ is real, thus for $y$ not the ordinate
of a zero, $\sin\theta = 0$ and $\cos\theta = \pm 1$.
Thus $\cos \theta $ is a discontinuous function. 
Now let $\yzero$ be the ordinate of a \emph{simple zero}. Then close
to such a zero we define 
\beq
c(y) =  \frac{ \chi (\tfrac{1}{2} + i y ) }{|\chi (\tfrac{1}{2} + i y)|}  
=  \frac{ y-\yzero}{|y - \yzero|}.
\eeq
For $y > \yzero$ then $c(y)=1$, and for $y<\yzero$ then  $c(y) = -1$,  
thus $c(y)$ is discontinuous precisely at a zero.
In the above polar representation, formally 
$c(y) =  \cos \theta (\tfrac{1}{2}, y)$. 
Therefore, by identifying zeros as the solutions to $\cos \theta =0$,  
we are simply defining  the function  $c(y)$ at the discontinuity 
as $c(\yzero)=0$. This is precisely what is displayed in 
\fref{fig:cos}, where the small $\delta$ smooths  out the discontinuity.
\end{remark}

\begin{remark}\label{general_sol}
On the critical line, but \emph{not} on a Riemann zero, 
$\theta = n\pi$ since $\chi=Ae^{i\theta}$ is real. Then
$\cos\theta = \pm 1$ and $\sin\theta = 0$. This shows
that $\cos\theta$ alternates in sign around a zero, being
a discontinuous function. The $\delta\to0^+$ limit smooths out this 
discontinuity so we can define $\cos\theta=0$ exactly on a zero $z=\rho$, 
where we also have $\sin\theta = \pm 1$.  
This can be confirmed numerically as illustrated in \fref{fig:cos} for the
first two zeros. 
\end{remark}

\begin{remark}\label{newzeta}
It is possible to introduce a new function
$\zeta(z)\mapsto\widetilde{\zeta}(z) = f(z) \zeta(z)$ that also 
satisfies the functional equation
\eqref{chisym}, i.e. $\tilde{\chi}(z)=\tilde{\chi}(1-z)$, but has zeros 
off the critical line due to the zeros of $f(z)$. In such a case 
the corresponding functional equation will hold if and only 
if $f(z) = f(1-z)$ for any $z$, and this is a trivial condition 
on $f(z)$, which could have been canceled in the first place. 
Moreover, if $f(z)$ and $\zeta(z)$ have different zeros, the analog of 
equation \eqref{sumchi} has a  factor
$f(z)$, i.e. $\tilde{\chi}(\rho+\delta) + \tilde{\chi}(1-\rho-\delta) = 
f(\rho+\delta)\[\chi(\rho+\delta)+\chi(1-\rho-\delta)\]=0$,  
implying \eqref{sumchi} again  where $\chi(z)$ is the original \eqref{chidef}. 
Therefore, the previous analysis eliminates $f(z)$
automatically and only finds the zeros of $\chi(z)$.
The analysis is  non-trivial precisely  because $\zeta(z)$ satisfies 
the functional equation but $\zeta(z) \ne \zeta(1-z)$. Furthermore, 
it is a well known theorem that the only function which satisfies the 
functional equation \eqref{chisym} and has the same characteristics 
of $\zeta(z)$, is $\zeta(z)$ itself. 
In other words, if $\widetilde{\zeta}(z)$ is required to have  the 
same properties of $\zeta(z)$, then $\widetilde{\zeta}(z) = C \, \zeta(z)$, 
where $C$ is a constant \cite[pg. 31]{Titchmarsh}.
\end{remark}

\begin{remark}\label{counting}
Although equations \eqref{exact_eq2} and \eqref{backlund} have an obvious 
resemblance, it is impossible to derive the former from the later, 
since the later is just a counting formula valid on the entire strip, 
and it is assumed that $T$ is \emph{not} the ordinate of a zero.
Moreover, this would require the assumption of the validity of the RH,
contrary to our approach, where we derived equations \eqref{exact_eq2}
and \eqref{FinalTranscendence} on the critical line, without assuming
the RH.
Despite our best efforts, we were not able to find 
formula \eqref{FinalTranscendence} in the literature. 
The formula \eqref{counting2} has never been proven on the critical 
line \cite{Edwards}. The current best estimate for the number of zeros on 
the critical line is given by $N_0(T) \ge \tfrac{2}{5} \, N(T)$ \cite{Conrey2}.
\end{remark}

\begin{remark} \label{circular}
One may object that our basic equation \eqref{FinalTranscendence} 
involves $\zeta(z)$ itself and this is somehow circular.
This is not a valid counter-argument.
First of all, $\arg \zeta$ already appears in the counting function $N(T)$.  
Secondly, the equation \eqref{FinalTranscendence} is a much more 
detailed equation than simply $\zeta (z )= 0$,  which has an infinite 
number of solutions,  in contrast with  \eqref{FinalTranscendence} which 
for each $n$ has a unique solution corresponding to the $n$-th zero.
Also, there are well-known ways to calculate the $\arg \zeta $ term,  
for example from an integral representation or
a convergent series \cite{Borwein}. 
\end{remark} 

\begin{figure}
\centering
\begin{minipage}{.5\textwidth}
  \centering
  \includegraphics[width=1\linewidth]{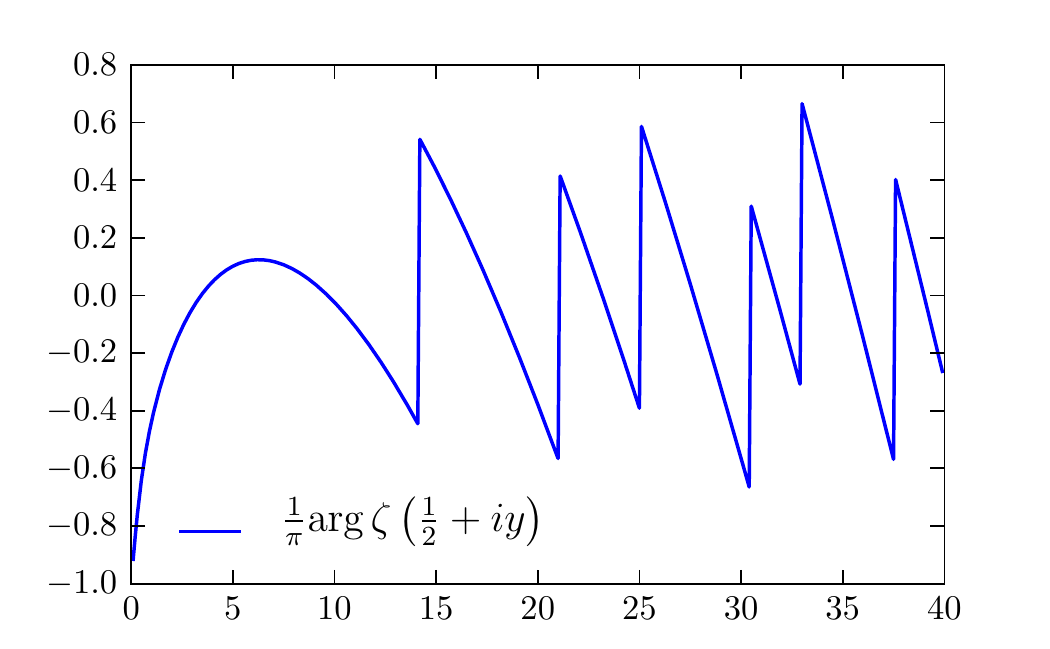}\\[-1em]
(a)
\end{minipage}%
\begin{minipage}{.5\textwidth}
  \centering
  \includegraphics[width=1\linewidth]{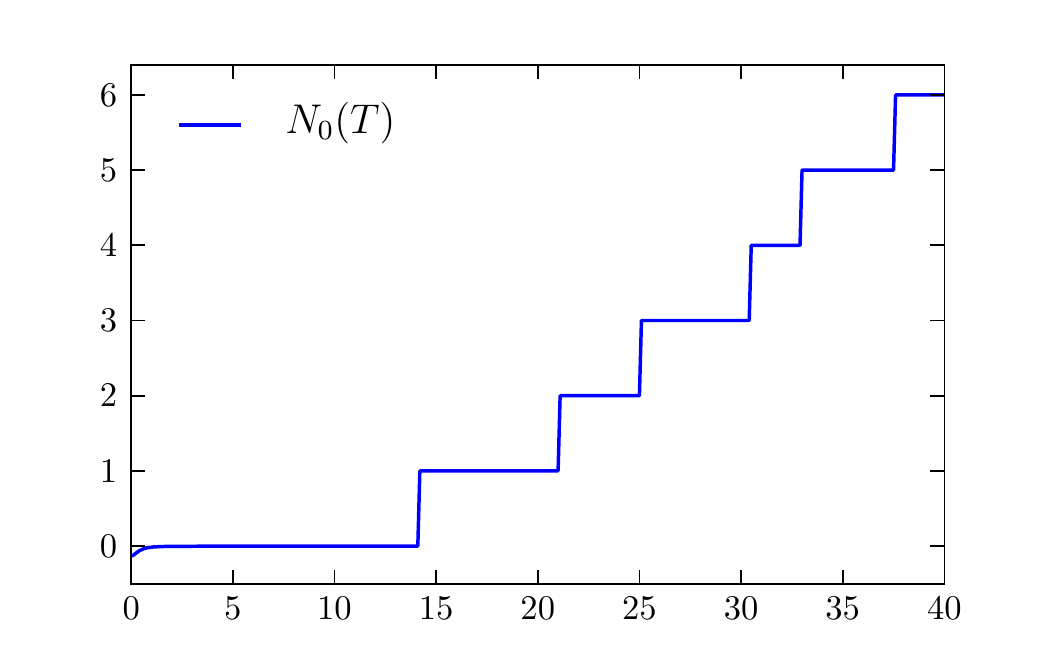}\\[-1em]
(b)
\end{minipage}
\caption{(a) A plot of $S(y)=\tfrac{1}{\pi}\arg \zeta\(\tfrac{1}{2} + i y\)$  
as a function of $y$ showing 
its rapid oscillation. The jumps occur on a Riemann zero. 
(b) The function $N_0(T)$ in \eqref{counting2}, which is indistinguishable 
from a manual counting of zeros.}
\label{fig:arg_counting}
\end{figure}

\begin{remark}
\label{arg_term}  
The small shift by $\delta$ in \eqref{FinalTranscendence} is essential
since it smooths out 
$S(y)  =  \tfrac{1}{\pi}\arg \zeta\(\tfrac{1}{2}+ i y\)$,  which is known 
to jump discontinuously at each zero. 
As well known, $\arg \zeta\(\tfrac{1}{2}+ i y\)$ 
is a piecewise continuous function, but rapidly oscillates around 
zero with discontinuous jumps, as shown in \fref{fig:arg_counting}a.
However, when this term  is added to the smooth part 
of $N(T)$, one obtains an accurate staircase function, which 
jumps by one at each zero on the line; see \fref{fig:arg_counting}b.  
 In this form, 
the formula \eqref{counting2} counts the zeros on the critical line 
accurately, i.e. it does not miss any zero.  Thus, as previously stated, 
since the Riemann-von Mangoldt function $N(T)$ has only been derived  
on the entire strip, and we have derived it for the zeros on the critical
line, this indicates that  all zeros are 
on the line.   
\end{remark}

\section{Approximate  solution in terms of the Lambert function}
\label{sec:lambert}

\subsection{Main formula} 

Let us now show that if one neglects the $\arg \zeta$ term,  
the equation \eqref{FinalTranscendence} can be exactly solved.
First, let us introduce the Lambert $W$ function \cite{Corless}, which is 
defined for any complex number $z$ through the equation
\beq
\label{lambert_definition}
W(z) e^{W(z)} = z.
\eeq
The multi-valued $W$ function cannot be expressed in terms of other
known elementary functions. 
If we restrict attention to real-valued $W(x)$ there are two branches. 
The principal branch occurs when $W(x) \ge -1$ and is denoted by $W_0$, 
or simply $W$ for short, and its domain is $x \ge -1/e$. 
The secondary branch, denoted by $W_{-1}$, satisfies $W_{-1}(x) \le -1$ 
for $-e^{-1} \le x < 0$.
Since we are interested in positive real-valued solutions  of 
\eqref{FinalTranscendence}, we just need the principal 
branch where $W$ is single-valued.

Let us consider the leading order approximation of 
\eqref{FinalTranscendence}, or equivalently, its average since
$\langle \arg\zeta\(\tfrac{1}{2}+iy\) \rangle = 0$. Then we have
the transcendental equation
\beq \label{ApproxTranscendence}
\dfrac{\tilde{y}_n}{2\pi}\log\(\dfrac{\tilde{y}_n}{2\pi e}\) = 
n - \dfrac{11}{8}.
\eeq
Through the transformation 
$\tilde{y}_n = 2\pi\(n-\tfrac{11}{8}\)x_n^{-1}$, this
equation can be written as 
$x_n e^{x_n} = e^{-1} \(n-\tfrac{11}{8}\)$. Comparing with 
\eqref{lambert_definition} its solution is  given by
$x_n = W\[e^{-1}\(n-\tfrac{11}{8}\)\]$, and thus we obtain 
\beq \label{Lambert}
\tilde{y}_n = 
\dfrac{2\pi\(n-\tfrac{11}{8}\)}{W\[e^{-1}\(n-\tfrac{11}{8}\)\]}.
\eeq

Although the inversion from \eqref{ApproxTranscendence} to \eqref{Lambert}
is rather simple, it is very  convenient since  it is indeed an 
explicit formula depending only on $n$,  
and $W$ is included in most numerical packages.
It gives an approximate solution for the ordinates of the Riemann zeros
in closed form. The values computed from \eqref{Lambert} are much closer 
to the Riemann zeros than Gram points,  
and one does not have to deal with violations of 
Gram's law (see below).     

\subsection{Further remarks}

\begin{remark}
The estimates given by \eqref{Lambert} can be calculated to high accuracy for 
arbitrarily large $n$, since $W$ is a standard elementary function.  
Of course, the  $\tilde{y}_n$ are not as accurate as the solutions 
$y_n$ including
the $\arg \zeta$ term, as we will see in section \ref{sec:numerical}.
Nevertheless, it is indeed a good estimate, especially if one considers 
very high zeros,  where  traditional methods have 
not previously estimated  such high values.
For instance, formula \eqref{Lambert} can easily estimate the 
 zeros shown in \tref{highn}, and much higher if desirable.
The numbers in this table are accurate approximations to the $n$-th zero 
to the number of digits shown,  which is approximately the number of 
digits in the integer part.
For instance,  the approximation to the $10^{100}$ zero is
correct to $100$ digits.   
With Mathematica we easily calculated the first million digits
of the $10^{10^6}$ zero.   
\end{remark}

\begin{table}
\centering
\begin{tabular}{@{}ll@{}}
\toprule[1pt]
$n$ & $\tilde{y}_n$ \\ 
\midrule[0.4pt]
$10^{22}+1$ &  
$1.370919909931995308226770\times 10^{21}$ \\
$10^{50}$ &  
$5.741532903784313725642221053588442131126693322343461\times 10^{48}$ \\
$10^{100}$  & 
$2.80690383842894069903195445838256400084548030162846045192360059224930$
\\[-1.6ex]
&$922349073043060335653109252473234\times 10^{98}$ \\
$10^{200}$ &
$1.38579222214678934084546680546715919012340245153870708183286835248393$
\\[-1.6ex]
&$8909689796343076797639408172610028651791994879400728026863298840958091$
\\[-1.6ex]
&$288304951600695814960962282888090054696215023267048447330585768
\times10^{198}$ \\
\bottomrule[1pt]
\end{tabular}
\caption{Formula \eqref{Lambert} can easily estimate very high Riemann zeros.
The results are expected to be correct up to the decimal point,  i.e. to the 
number of digits in the integer part.   
The numbers are shown with three digits beyond the integer part.}
\label{highn}
\end{table}

\begin{remark}Using the asymptotic behaviour $W(x) \sim \log x$ for large $x$,
the $n$-th zero is approximately $\tilde{y}_n \approx  2 \pi n/ \log n$,
as already known \cite{Titchmarsh}.
The distance between consecutive zeros is $2\pi/\log n$, which
tends to zero when $n\to \infty$.
\end{remark}

\begin{remark}\label{gram}
The solutions $\tilde{y}_n$ to the equation \eqref{ApproxTranscendence}
are reminiscent of the so-called Gram 
points $g_n$,  which are solutions to 
$\vartheta (g_n) = n\pi$ where $\vartheta$ is given by \eqref{riemann_siegel}. 
Gram's law is the tendency for Riemann zeros to lie between 
consecutive Gram points, but it is known to fail for about $1/4$ 
of all 
Gram intervals. Our $\tilde{y}_n$ are intrinsically different from Gram 
points, being an approximation to the ordinate of the $n$-th
Riemann zero. In particular, the  Gram point $g_0 = 17.8455$ is the closest to 
the first Riemann zero, whereas, $\tilde{y}_1 = 14.52$ is 
much closer to the true zero which is $y_1 \simeq 14.1347$.
The traditional method to compute the zeros is based on 
the Riemann-Siegel formula, 
$\zeta\(\tfrac{1}{2}+it\) = Z(t)\(\cos\vartheta(t) - i \sin\vartheta(t)\)$,
and the empirical observation that the real part of this equation 
is almost always positive, except when Gram's law fails, and $Z(t)$ has 
the opposite
sign of $\sin\vartheta$. Since $Z(t)$ and $\zeta\(\tfrac{1}{2}+it\)$ have
the same zeros, one looks for the zeros of $Z(t)$ between two
Gram points, as long as Gram's law holds $(-1)^nZ\(g_n\)>0$. 
To verify the RH numerically, the counting formula \eqref{backlund} must 
also be used, to assure that the number of zeros on the critical line 
coincide with the number of zeros on the strip. 
The detailed procedure is throughly 
explained in \cite{Edwards,Titchmarsh}.
Based on this method, amazingly accurate solutions and high zeros on 
the critical line were computed 
\cite{OdlyzkoSchonhage,Odlyzko,Odlyzko2,Gourdon}.
Nevertheless, our proposal is \emph{fundamentally} different.
We claim that \eqref{exact_eq2}, or its asymptotic 
approximation \eqref{FinalTranscendence}, is the equation that determines
the Riemann zeros on the critical line.
Then, one just needs to find its solution for a given
$n$. We will compute the Riemann zeros  in this way in the next 
section, just by solving the equation numerically, starting
from the approximation given by the explicit formula \eqref{Lambert},
without using Gram points nor the Riemann-Siegel $Z$ function.
Let us emphasize that our goal is not to provide a more efficient algorithm
to compute the zeros \cite{OdlyzkoSchonhage},  although the method 
described here may very well be,   but to justify the 
validity of equations \eqref{FinalTranscendence} and \eqref{exact_eq2}.
\end{remark}

\section{Numerical solutions}
\label{sec:numerical}

Instead of solving the exact equation \eqref{exact_eq2} we will 
initially consider its first order approximation, 
which is equation \eqref{FinalTranscendence}. 
As we will see, this approximation already yields surprisingly accurate values
for the Riemann zeros.

\begin{figure}
\centering
\begin{minipage}{.5\textwidth}
  \centering
  \includegraphics[width=1\linewidth]{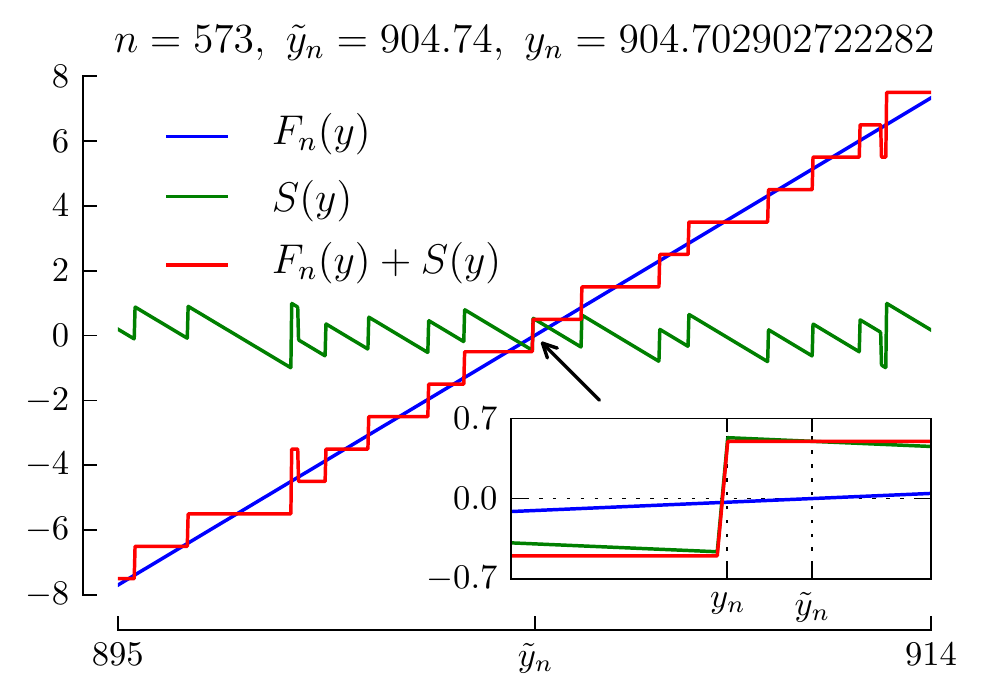}\\[-1em]
(a)
\end{minipage}%
\begin{minipage}{.5\textwidth}
  \centering
  \vspace{1.1em}
  \includegraphics[width=1\linewidth]{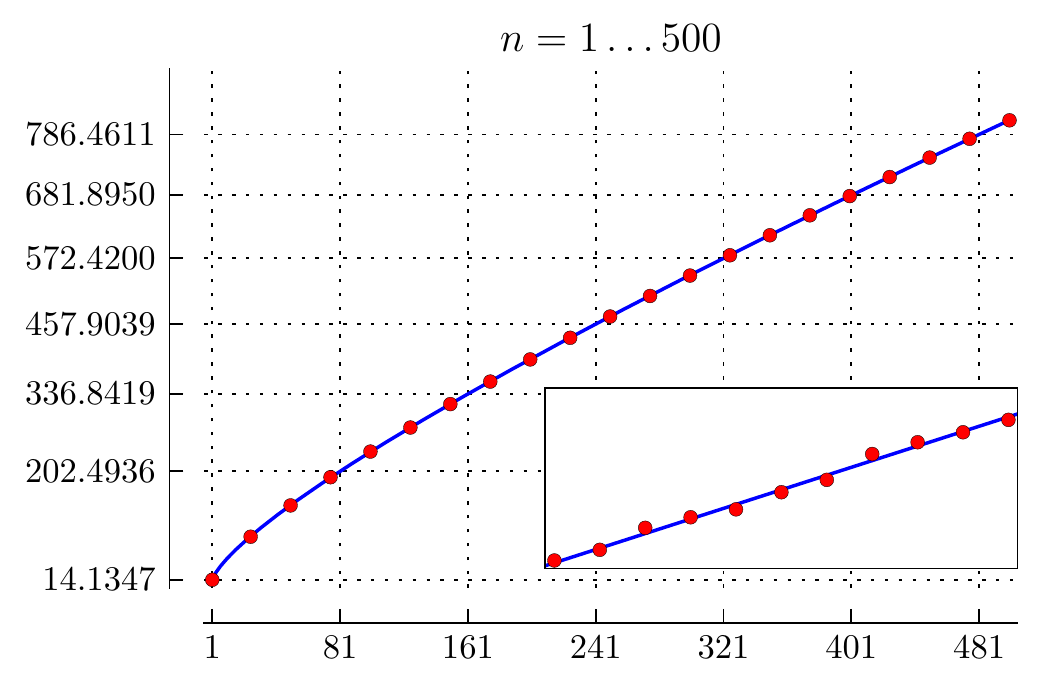}\\[-1em]
(b)
\end{minipage}
\caption{
(a) Graph of the different terms present in equation
\eqref{FinalTranscendence}. Here we have
$F_n(y) = \tfrac{y}{2\pi}\log\(\tfrac{y}{2\pi e}\) - n +\tfrac{11}{8}$, 
equation \eqref{ApproxTranscendence}, and 
$S(y) = \tfrac{1}{\pi}\arg\zeta\(\tfrac{1}{2}+iy\)$. $F_n(y)+S(y)$ is the
complete equation \eqref{FinalTranscendence}.
Note how $\tilde{y}_n$ is dislocated by 
the $\arg \zeta$ term, yielding a much more precise result.
(b) Comparison of the prediction of \eqref{Lambert} (blue line) and 
\eqref{FinalTranscendence} (red dots). The numerical solutions of 
\eqref{FinalTranscendence} oscillate around the line predicted by 
\eqref{Lambert} due to the fluctuating term $\arg \zeta$.
}
\label{fig:trans_zeros}
\end{figure}

Let us first consider how the approximate solution given by
\eqref{Lambert} is modified by the presence of the $\arg\zeta$ term
in \eqref{FinalTranscendence}.  Numerically, we compute $\arg\zeta$ taking
its principal value.
As already discussed in Remark \ref{arg_term}, the function 
$\arg\zeta\(\tfrac{1}{2} + i y\)$ oscillates around zero and changes sign 
in the vicinity of each Riemann zero, as shown in \fref{fig:arg_counting}a.
At a zero it can be well-defined by the limit \eqref{deltadef}, 
which is generally not zero. For example, for the first Riemann zero 
$y_1\simeq 14.1347$, 
$\lim_{\delta \to 0^+} \arg \zeta \( \tfrac{1}{2}+\delta +  i y_1 \) = 
0.1578739$.
The $\arg \zeta$ term plays an important role and 
indeed improves the estimate of the $n$-th zero. This can be seen from
\fref{fig:trans_zeros}a for a randomly chosen $n$. It  
practically cancels \eqref{ApproxTranscendence} around
the zero, and exactly at the true zero we have a jump.
The value predicted by \eqref{Lambert} is  then slightly changed. 
For a given $n$, the problem of finding the value $y_n$ 
where this jump occurs, yields the $n$-th Riemann 
zero as the numerical solution of \eqref{FinalTranscendence}.

Since equation \eqref{FinalTranscendence} alternates in sign 
around a zero, it is convenient to use Brent's method \cite{BrentMethod} 
to find its root. We applied this method, looking for a root in an appropriate
interval, centered around the approximate solution $\tilde{y}_n$ given 
by formula \eqref{Lambert}.  Some of the solutions  are presented
in \tref{some_zeros}, and are accurate up to the number of decimal places
shown. We used  only Mathematica or some very simple algorithms to perform 
these  numerical computations,
taken from standard open source numerical libraries.
We present the numbers accurate up to $9$ digits after the 
integer part. 

\begin{table}
\centering
\begin{tabular}{@{}lrr@{}}
\toprule[1pt]
$n$ & $\tilde{y}_n$ & $y_n$ \\ 
\midrule[0.4pt] 
$1$       &         $14.52$ &         $14.134725142$ \\
$10$      &         $50.23$ &         $49.773832478$ \\
$10^2$    &        $235.99$ &        $236.524229666$ \\
$10^3$    &       $1419.52$ &       $1419.422480946$ \\
$10^4$    &       $9877.63$ &       $9877.782654006$ \\
$10^5$    &      $74920.89$ &      $74920.827498994$ \\
$10^6$    &     $600269.64$ &     $600269.677012445$ \\
$10^7$    &    $4992381.11$ &    $4992381.014003180$ \\
$10^8$    &   $42653549.77$ &   $42653549.760951554$ \\
$10^9$    &  $371870204.05$ &  $371870203.837028053$ \\
$10^{10}$ & $3293531632.26$ & $3293531632.397136704$ \\
\bottomrule[1pt]
\end{tabular}
\caption{Numerical solutions of equation \eqref{FinalTranscendence}.
All numbers shown are accurate up to the $9$-th decimal place, in 
comparison with  \cite{Odlyzko,Oliveira}.}   
\label{some_zeros}
\end{table}

\begin{table}
\centering
\begin{tabular}{@{}ll@{}}
\toprule[1pt]
$n$ & $y_n$ \\ 
\midrule[0.4pt] 
$1$ & $14.13472514173469379045725198356247$ \\
$2$ & $21.02203963877155499262847959389690$ \\
$3$ & $25.01085758014568876321379099256282$ \\
$4$ & $30.42487612585951321031189753058409$ \\
$5$ & $32.93506158773918969066236896407490$ \\
\bottomrule[1pt]
\end{tabular}
\caption{Numerical solutions to  \eqref{FinalTranscendence}  for  the 
lowest zeros.   Although it  was derived for high $y$,
it provides accurate numbers even for the lower zeros. 
This numbers are correct up to the decimal place shown \cite{Odlyzko}.}
\label{lower_zeros}
\end{table}

Although the formula for $y_n$ was derived for large $y$, it is 
surprisingly accurate even for the lower zeros, as shown in
\tref{lower_zeros}.   It is actually easier to solve numerically for low 
zeros since  $\arg \zeta $ is better behaved.    These numbers are correct 
up to the number of digits shown, and the precision was improved simply 
by decreasing the error tolerance. 

Riemann zeros have previously been calculated to high accuracy using
sophisticated algorithms \cite{OdlyzkoSchonhage}, which are not based on 
solving our equation \eqref{FinalTranscendence}.        
Nevertheless,  we have verified that \eqref{FinalTranscendence} 
is well satisfied  to the degree of accuracy of these zeros. 
This can be seen in \tref{checking} where we show
the absolute value of \eqref{FinalTranscendence}, replaced with our 
numerical solutions,
and its value calculated with much more accurate Riemann zeros, 
up to the $150$-th decimal place, provided by Mathematica.

\begin{table}
\centering
\begin{minipage}{.49\textwidth}
\centering
\begin{tabular}{@{}lll@{}}
\toprule[1pt]
$n$ & $|\mbox{\eqref{FinalTranscendence}}|$ & 
$|\mbox{\eqref{FinalTranscendence}}|$ \\
\midrule[0.4pt] 
$1$ & $0.082$ & $0.00047$ \\
$2$ & $0.071$ & $0.00032$ \\
$3$ & $0.046$ & $0.00026$ \\
$4$ & $0.044$ & $0.00022$ \\
$5$ & $0.081$ & $0.00020$ \\
\bottomrule[1pt]
\end{tabular}
\end{minipage}
\begin{minipage}{0.49\textwidth}
\centering
\begin{tabular}{@{}lll@{}}
\toprule[1pt]
$n$ & $|\mbox{\eqref{FinalTranscendence}}|$ & 
$|\mbox{\eqref{FinalTranscendence}}|$ \\
\midrule[0.4pt] 
$10^5-4$ & $0.10$ & $0.000000088$ \\
$10^5-3$ & $0.094$ & $0.000000089$ \\
$10^5-2$ & $0.072$ & $0.000000090$ \\
$10^5-1$ & $0.034$ & $0.000000088$ \\
$10^5$ & $0.059$ & $0.000000089$ \\
\bottomrule[1pt]
\end{tabular}
\end{minipage}
\caption{
We substitute numerical solutions of $y_n$ into equation
\eqref{FinalTranscendence}.  The values $|(\ref{FinalTranscendence})|$  
denote the absolute value of the  difference
between the RHS and LHS of the equation.    
We use $\delta=10^{-9}$. The second column is the value obtained with our 
solutions, accurate to 9 decimal places. The third column is the value 
obtained with solutions accurate up to 150 decimal places, computed through 
another algorithm, which is  not based on solving \eqref{FinalTranscendence}. 
This shows that \eqref{FinalTranscendence} is  indeed  satisfied for 
high accurate Riemann zeros.
}
\label{checking}
\end{table}

\section{GUE Statistics}
\label{sec:gue}

The link between the Riemann zeros and random matrix theory
started with the pair correlation of zeros, proposed by 
Montgomery \cite{Montgomery},  and
the observation of F. Dyson that it is the same as the 2-point correlation
function predicted by the gaussian unitary ensemble (GUE) for large 
random  matrices \cite{Dyson}.

The main purpose of this section is to test whether our 
approximation \eqref{FinalTranscendence} to the zeros 
is accurate enough to reveal this statistics.
Whereas formula \eqref{Lambert} 
is a valid estimate of  the zeros, 
it is not sufficiently accurate to reproduce 
the GUE statistics,  since it does not have the 
oscillatory  $\arg \zeta$ term. 
On the other hand, the solutions to equation \eqref{FinalTranscendence} 
are accurate enough,
which indicates the importance of the $\arg \zeta$.

Montgomery's pair correlation conjecture can be stated as follows:
\beq
\label{montgomery}
\dfrac{1}{N(T)}
\sum_{\substack{
0\le y,y'\le T \\[0.4em]
\alpha < \tfrac{1}{2\pi}\log\(\tfrac{T}{2\pi}\)\(y-y'\) \le \beta
}} 1
\quad \sim \quad  \int_{\alpha}^{\beta}du
\(1 - \dfrac{\sin^2\(\pi u\)}{\pi^2 u^2}\)
\eeq
where $ 0 < \alpha<\beta$, 
$N(T)\sim \tfrac{T}{2\pi}\log\(\tfrac{T}{2\pi}\)$ 
according to \eqref{counting2}, and the statement is valid in the limit
$T\to \infty$. The
right hand side of \eqref{montgomery} is the 2-point GUE correlation 
function. The average spacing between consecutive zeros is given by 
$\tfrac{T}{N} \sim 2\pi/\log\(\tfrac{T}{2\pi}\)\to 0$  
as $T\to \infty$. This can also be seen from \eqref{Lambert} for very
large $n$, i.e. $\tilde{y}_{n+1}-\tilde{y}_n \to 0$ as $n\to\infty$. 
Thus the distance between zeros on the
left hand side of \eqref{montgomery}, under the sum, is a normalized
distance.

While \eqref{montgomery} can be applied if we start from the first
zero on the critical line, it is unable to provide a test if we are centered
around a given high zero on the line. To deal with such a situation,
Odlyzko \cite{Odlyzko2} proposed a stronger version of Montgomery's 
conjecture, by taking into account the large density of zeros 
higher on the line. This is done by replacing the normalized distance in 
\eqref{montgomery} by a sum of normalized distances over
consecutive zeros in the form
\beq
\delta_n \equiv \dfrac{1}{2\pi}\log\(\dfrac{y_n}{2\pi}\)\(y_{n+1}-y_n\).
\eeq
Thus \eqref{montgomery} is replaced by
\beq
\label{odlyzko_pair}
\dfrac{1}{\(N-M\)\(\beta-\alpha\)}\sum_{\substack{
M \le m,n \le N \\[0.4em]
\alpha < \sum_{k=1}^{n}\delta_{m+k} \le \beta
}} 1
\quad = \quad  \dfrac{1}{\beta-\alpha}\int_{\alpha}^{\beta}du
\(1 - \dfrac{\sin^2\(\pi u\)}{\pi^2 u^2}\),
\eeq
where $M$ is the label of a given zero on the line and
$N > M$. In this sum it is assumed that $n > m$ also, and we
included the correct normalization on both sides. The conjecture
\eqref{odlyzko_pair} is already  well supported by extensive 
numerical analysis \cite{Odlyzko2,Gourdon}. 

\begin{figure}
\centering
\begin{minipage}{.49\textwidth}
  \centering
  \includegraphics[width=1\linewidth]{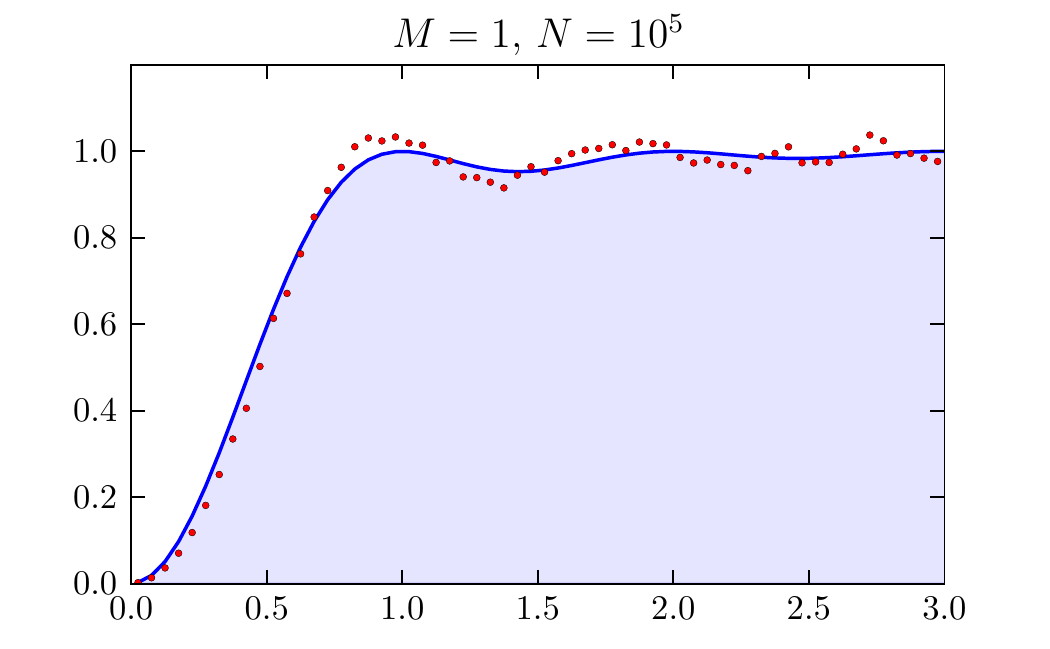}\\[-1em]
(a)
\end{minipage}
\begin{minipage}{.49\textwidth}
  \centering
  \includegraphics[width=1\linewidth]{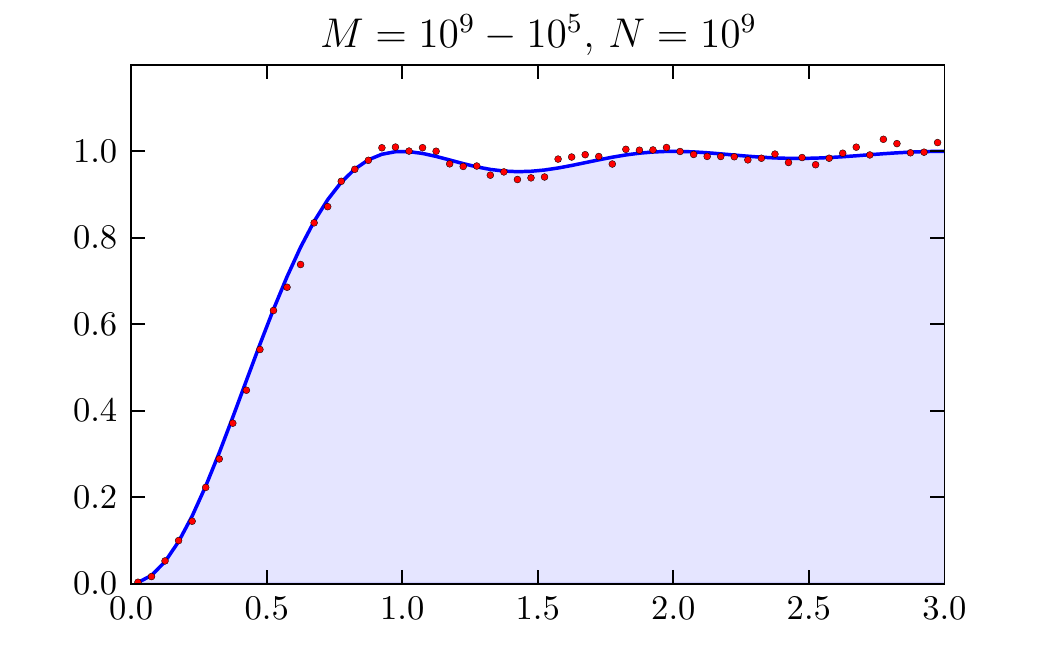}\\[-1em]
(b)
\end{minipage}
\caption{
(a) The solid line represents the RHS of \eqref{odlyzko_pair}
and the dots represent its LHS, computed from the
numerical solutions of equation \eqref{FinalTranscendence}.
The parameters are $\beta = \alpha + 0.05$,
$\alpha=[0,\, 0.05,\,\dotsc,\,3]$ and the $x$-axis is given by
$x=\tfrac{1}{2}\(\alpha+\beta\)$. We use the first $10^5$ zeros.
(b) The same parameters but using zeros in the middle of the
critical line, i.e. $M=10^9-10^5$ and $N=10^9$.
}
\label{fig:gue}
\end{figure}

Odlyzko's conjecture \eqref{odlyzko_pair} is a very strong constraint on
the statistics of the zeros. Thus we submit the numerical solutions of 
equation \eqref{FinalTranscendence}, as discussed
in the previous section, to this test. In 
\fref{fig:gue}a we can see the result for $M=1$ and $N=10^{5}$, with
$\alpha$ ranging from $0\dotsc 3$ in steps of $s=0.05$, and
$\beta=\alpha+s$ for each value of $\alpha$, i.e. 
$\alpha = [0.00, \, 0.05, \, 0.10,  \dotsc,\,3.00]$ and 
$\beta = [0.05,\, 0.10, \, \dotsc, \, 3.05]$. 
We compute the left hand side of \eqref{odlyzko_pair} for each 
pair $(\alpha, \beta)$ and plot
the result against $x = \tfrac{1}{2}\(\alpha + \beta\)$.
In \fref{fig:gue}b we do the same thing but with
$M=10^9-10^5$ and $N=10^9$.
Clearly, the numerical solutions of \eqref{FinalTranscendence} reproduce
the correct statistics. In fact, \fref{fig:gue}a is identical
to the one in \cite{Odlyzko2}. The last zeros in these ranges are
shown in \tref{high_values}.

\begin{table}
\centering
\begin{minipage}{.49\textwidth}
\centering
\begin{tabular}{@{}ll@{}}
\toprule[1pt]
$n$ & $y_n$ \\ 
\midrule[0.4pt] 
$10^5-5$ & $74917.719415828$ \\
$10^5-4$ & $74918.370580227$ \\
$10^5-3$ & $74918.691433454$ \\
$10^5-2$ & $74919.075161121$ \\
$10^5-1$ & $74920.259793259$ \\
$10^5$   & $74920.827498994$ \\
\bottomrule[1pt]
\end{tabular}
\end{minipage}
\begin{minipage}{.49\textwidth}
\centering
\begin{tabular}{@{}ll@{}}
\toprule[1pt]
$n$ & $y_n$ \\ 
\midrule[0.4pt] 
$10^9-5$ & $371870202.244870424$ \\
$10^9-4$ & $371870202.673284411$ \\
$10^9-3$ & $371870203.177729726$ \\
$10^9-2$ & $371870203.274345875$ \\
$10^9-1$ & $371870203.802552402$ \\
$10^9$   & $371870203.837028146$ \\
\bottomrule[1pt]
\end{tabular}
\end{minipage}
\caption{Last numerical solutions to \eqref{FinalTranscendence} around
$n=10^5$ and $n=10^9$. In the first table the solutions are accurate 
up to the $8$-th decimal place, while in the second table up to the 
$6$-th decimal place \cite{Odlyzko,Oliveira}.}
\label{high_values}
\end{table}

\section{Prime number counting function}  
\label{sec:prime}

In this section we explore whether our approximations to the Riemann zeros
are accurate enough to reconstruct the prime number counting function.     
As usual, let $\pi (x)$ denote the number of primes less than $x$.
Riemann obtained an explicit expression for $\pi (x)$ in terms of 
the non-trivial zeros of $\zeta(z)$.  
There are simpler but equivalent versions 
of the main result,  based  on the function $\psi (x) $ below.   
However, let us present the main formula for $\pi (x)$ itself,
since it is historically more important.

The function $\pi (x)$ is related to another number-theoretic 
function $J(x)$, defined as 
\beq
\label{Jx}
J(x)  =  \sum_{2\leq n \leq x}    \frac{\Lambda (n)}{\log n} 
\eeq
where  $\Lambda (n)$,  the  von Mangoldt function,
is equal to $\log p$ if $n=p^m$ for some prime $p$  
and an integer $m$, and zero otherwise.  The two functions 
$\pi (x)$ and $J(x)$ are related by M\" obius inversion:
\beq
\label{mob1}
\pi (x) = \sum_{n\geq 1} \frac{\mu (n)}{n}  J(x^{1/n}).
\eeq
Here,  $\mu (n)$ is the M\" obius function,  
equal to $1$ ($-1$)  if $n$ is a product of an even (odd) number
of distinct primes,  
and equal to zero if it has a multiple prime factor.    
The above expression is actually a finite sum,  
since for large enough $n$,  $x^{1/n} <2$ and $J=0$.

The main result of Riemann is a formula for $J(x)$, expressed 
as an infinite sum over zeros $\rho$  of the $\zeta(z)$ function:
\beq
\label{Jzeros}
J(x) =  \Li (x) - \sum_\rho  \Li\(x^\rho\)  +  
\int_x^\infty  \dfrac{dt}{\log t} ~  \inv{t \(t^2 -1\)} - \log 2,
\eeq
where $\Li (x) = \int_0^x dt /\log t $ is the log-integral 
function \footnote{Some care must be taken in  numerically 
evaluating $\Li (x^\rho )$  since $\Li$ has a branch point. 
It is more properly defined as ${\rm Ei} (\rho \log x )$ 
where ${\rm Ei} (z) = - \int_{-z}^\infty  dt \, e^{-t} /t$ is the 
exponential integral function.}.
The above sum is real because the $\rho$'s come 
in conjugate pairs.
If there are no zeros on the line $\Re(z) = 1$, then the dominant 
term is the first one in the above equation, $J(x) \sim  \Li (x)$,  
and this was used to  prove the prime number theorem by 
Hadamard and de la Vall\' ee Poussin.

The function $\psi (x)$ has the simpler form
\beq
\label{psizeros}
\psi (x)  =   \sum_{n \leq x}   \Lambda (n) = 
x - \sum_\rho\dfrac{x^\rho}{\rho} - \log (2\pi)  - \inv{2} 
\log \( 1 - \inv{x^2} \).
\eeq
In this formulation,  the prime number theorem follows from the fact 
that the leading term is $\psi (x) \sim x$.

\begin{figure}
\centering
\begin{minipage}{.5\textwidth}
  \centering
  \includegraphics[width=1\linewidth]{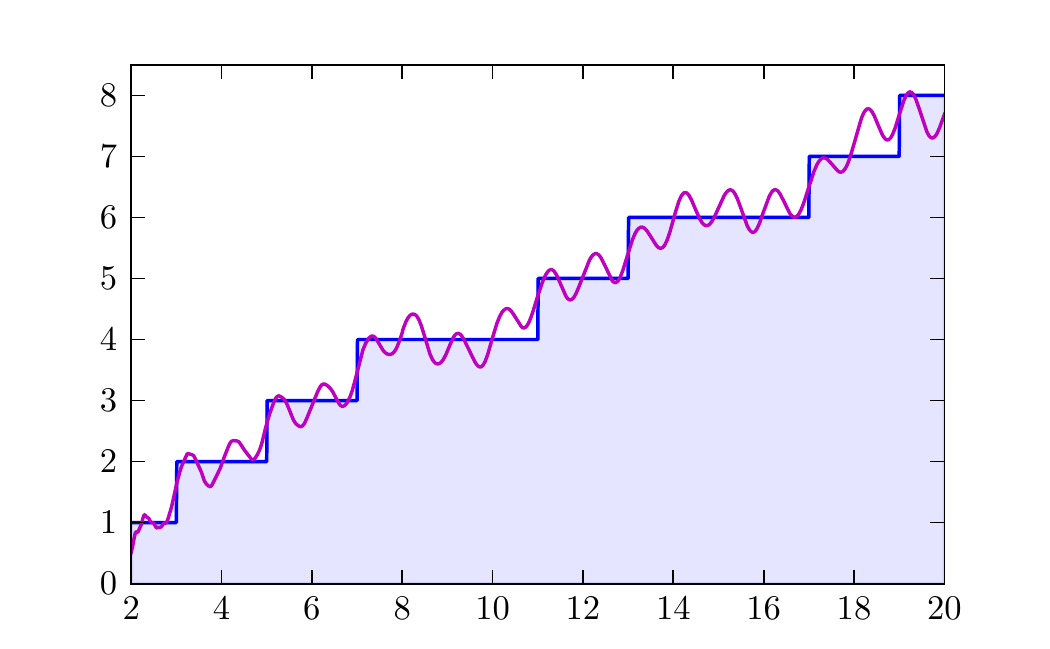}\\[-1em]
(a)
\end{minipage}%
\begin{minipage}{.5\textwidth}
  \centering
  \includegraphics[width=1\linewidth]{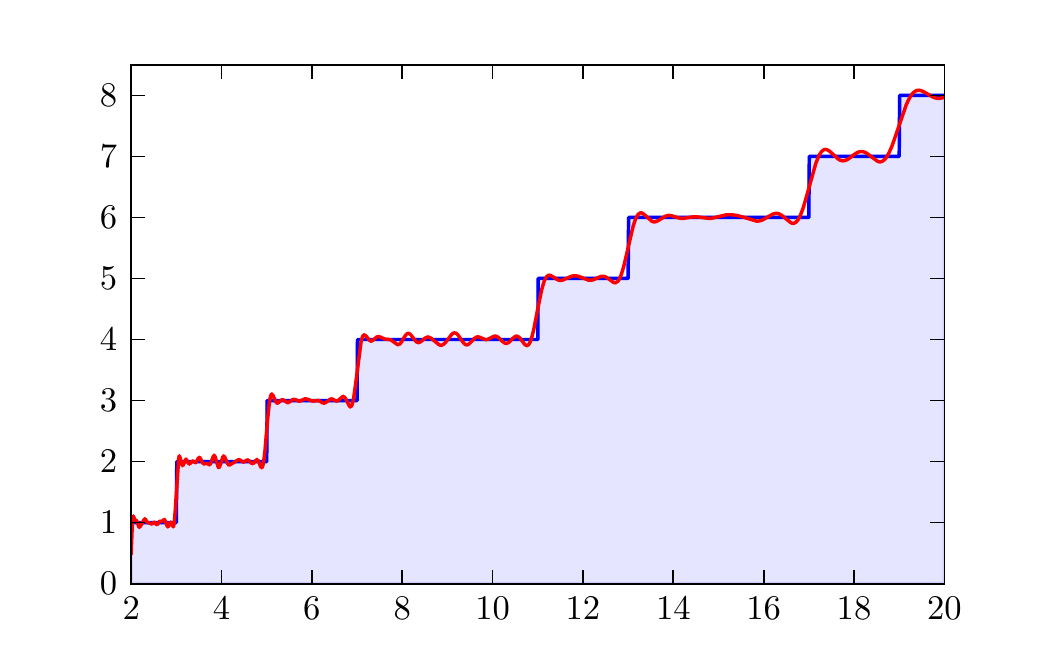}\\[-1em]
(b)
\end{minipage}
\caption{
(a) 
The prime number counting function $\pi (x)$ with the first $50$ Riemann
zeros approximated by the formula \eqref{Lambert}.
(b) 
The same plot but with the first $50$ Riemann zeros approximated 
by numerical solutions to the equation \eqref{FinalTranscendence}. 
}
\label{fig:prime}
\end{figure}

In Figure \fref{fig:prime}a we plot $\pi(x)$  from equations \eqref{mob1}
and \eqref{Jzeros}, computed with the first $50$ zeros in the 
approximation $\rho_n = \tfrac{1}{2} + i \tilde{y}_n$ given by 
\eqref{Lambert}. \fref{fig:prime}b  shows the same plot with zeros 
obtained from the numerical solution of equation \eqref{FinalTranscendence}.
Although with the approximation $\tilde{y}_n$ the curve is trying 
to follow the steps in
$\pi (x)$,  once again, one clearly sees the importance 
of the $\arg  \zeta$ term.

\section{Numerical solutions to the exact equation}
\label{sec:numerical_exact}

In the previous sections we have computed numerical solutions of
\eqref{FinalTranscendence} showing that, actually, this first
order approximation to \eqref{exact_eq2} is very good and already 
captures  the  interesting properties of the Riemann zeros, 
like the GUE statistics and ability
to reproduce the prime number counting formula. Nevertheless, by simply 
solving \eqref{exact_eq2} it is possible to obtain values for the zeros
as accurately  as desirable. 
The numerical procedure is performed  as follows: 
\begin{enumerate}
\item \label{step1} We solve \eqref{exact_eq2} looking for the solution in 
a region centered around the number $\tilde{y}_n$ provided 
by \eqref{Lambert}, with a not so small $\delta$, 
for instance $\delta \sim 10^{-5}$.
\item \label{step2} We solve \eqref{exact_eq2} again but now centered around 
the solution obtained in step \ref{step1} above, and we decrease $\delta$, for 
instance $\delta \sim 10^{-8}$.
\item We repeat the procedure in step \ref{step2} above, 
decreasing $\delta$ again.
\item Through successive iterations, and decreasing $\delta$ 
each time, it is possible to obtain solutions as accurate as desirable.
In carrying this out,  
it is important to not allow $\delta$ to be exactly zero.   
\end{enumerate}

The first few zeros are shown in \tref{lower_precise}.
We simply applied the standard root finder  in 
Mathematica \footnote{The Mathematica notebook 
we used to carry out these computations has only a few dozen lines of code
and is available on the arXiv in math.NT as an auxiliary file to this 
submission.}. 
Through successive iterations it is possible achieve even much higher  
accuracy than  shown in 
\tref{lower_precise}.

\begin{table}
\centering
\begin{tabular}{@{}ll@{}}
\toprule[1pt]
$n$ & $y_n$ \\ 
\midrule[0.4pt] 
$1$ & $14.1347251417346937904572519835624702707842571156992431756855$ \\
$2$ & $21.0220396387715549926284795938969027773343405249027817546295$ \\
$3$ & $25.0108575801456887632137909925628218186595496725579966724965$ \\
$4$ & $30.4248761258595132103118975305840913201815600237154401809621$ \\
$5$ & $32.9350615877391896906623689640749034888127156035170390092800$ \\
\bottomrule[1pt]
\end{tabular}
\caption{The first few numerical solutions to \eqref{exact_eq2},
accurate to $60$ digits ($58$ decimals). These  solutions
were obtained using  the root finder function in Mathematica.}
\label{lower_precise}
\end{table}

It is known that the first zero where  Gram's law fails is for $n=126$.
Applying the same method, like for any other $n$, the solution of 
\eqref{exact_eq2} starting with the approximation \eqref{Lambert} does
not present any difficulty. We  easily found the following number:
\begin{flalign*}
279.229250927745189228409880451955359283492637405561293594727
\qquad (n=126)
\end{flalign*}
Just to illustrate, and to convince the reader,  how the solutions 
of \eqref{exact_eq2} can be made arbitrarily precise, we compute the 
zero $n=1000$ accurate up to $500$ decimal places, also using the same 
simple approach \footnote{Computing this number to $500$ digit accuracy  
took a few minutes on a standard 8 GB RAM laptop using
Mathematica. It only takes a few seconds to obtain 100 digit accuracy.}:
\begin{flalign*}
1419.&42248094599568646598903807991681923210060106416601630469081468460\\[-8pt]
&8676417593010417911343291179209987480984232260560118741397447952650637\\[-8pt]
&0672508342889831518454476882525931159442394251954846877081639462563323\\[-8pt]
&8145779152841855934315118793290577642799801273605240944611733704181896\\[-8pt]
&2494747459675690479839876840142804973590017354741319116293486589463954\\[-8pt]
&5423132081056990198071939175430299848814901931936718231264204272763589\\[-8pt]
&1148784832999646735616085843651542517182417956641495352443292193649483\\[-8pt]
&857772253460088
\end{flalign*}
Substituting precise Riemann zeros calculated by other 
means \cite{Odlyzko} into \eqref{exact_eq2} one
can check that the equation is identically satisfied.
These results corroborate that \eqref{exact_eq2} is an
exact equation for the Riemann zeros, which was derived on the 
critical line.

\section{Final remarks}
\label{sec:conclusion}

Let us summarize our main results and arguments.
Throughout this paper we did \emph{not} assume the Riemann hypothesis.
The main outcome was the demonstration that there are infinite zeros
on the critical line, $\rho_n = \tfrac{1}{2} + i y_n$, where $y_n$ 
exactly satisfies the equation \eqref{exact_eq2}. Asymptotically this
equation can be approximated by \eqref{FinalTranscendence}.
Furthermore, we argued that these equations can be made continuous through
the $\delta\to 0^+$ limit, and therefore, they should have a unique solution for 
every single $n$. Under this assumption, the number of solutions
on the critical line already saturates the counting formula for the number
of zeros on the entire critical strip. This is a strong indication that
\eqref{exact_eq2} captures all non-trivial zeros, which must 
therefore be all on the critical line. Although our approach cannot be
considered as a rigorous proof,  it is at the very least a clear  
strategy  towards 
proving the Riemann hypothesis. 
It is important to note that \eqref{exact_eq2} and
\eqref{FinalTranscendence} were \emph{derived} on the critical line, while
the counting formulas \eqref{backlund} and \eqref{riemann_counting} can
only be derived on the entire strip. Thus it is impossible to obtain the 
former from the latter without assuming the Riemann hypothesis.

We verified numerically that the simplest approximation to 
the exact equation \eqref{exact_eq2}, namely \eqref{FinalTranscendence}, 
is enough to capture the statistical properties of the Riemann zeros. 
We did so by testing the Montgomery-Odlyzko pair correlation conjecture, 
and by reconstructing the prime number counting function, employing 
the numerical solutions of equation \eqref{FinalTranscendence}. 
In solving such transcendental equation, we started from an approximate 
solution given by the explicit formula \eqref{Lambert}. Thus, we did not
require the use of Gram points and we also did not have to deal with
violations of Gram's law. We also computed some numerical solutions
of the exact equation \eqref{exact_eq2}, which proved to be much more
stable under the numerical approach. This procedure constitutes a novel
method to compute the zeros.
Therefore, the numerical results strongly support the validity of our 
assertions, claiming  that \eqref{exact_eq2} is an exact equation, 
identically satisfied by the $n$-th Riemann zero on the critical line.

We also wish to mention that we have extended this work to two 
infinite classes of $L$-functions, those based on Dirichlet 
characters and modular forms \cite{Lfunctions}.

\section*{Acknowledgments}

AL wishes to thank Giuseppe Mussardo and Germ\' an Sierra for discussions.
We also  wish to  thank  Tim Healey,  Christopher Hughes,    
Wladyslaw Narkiewicz, 
Andrew Odlyzko,   Mark Srednicki,   and Tao Su
for critical comments on the first draft. AL  is grateful to 
the hospitality of the Centro  Brasileiro de Pesquisas  F\' isicas in Rio 
de Janeiro where this work was completed, especially Itzhak Roditi,
and the support of CNPq  under the ``Ci\^ encias sem fronteiras''
program in Brazil,  which also supports GF.  
This work is supported by the National Science 
Foundation of the United States of America under grant number NSF-PHY-0757868.

\end{document}